\documentclass[leqno]{amsart}
\usepackage{bbm}%***for using the black board meta fonts like upper half plane
\usepackage{amsmath,amsthm,amsfonts,amssymb}
\usepackage{color}
\usepackage{enumerate}

\newtheorem*{acknowledgements}{Acknowledgements}

\newtheorem*{remark0}{Remark}
\newtheorem*{theorem1}{Theorem 1}

\newtheorem*{note1}{Note}
\newtheorem*{Breteche}{Remark}

\newtheorem*{theorem2}{Theorem 2}
\newtheorem*{corollary3}{Corollary 3}

\newtheorem*{Conjecture}{Mass Equidistribution}
\newtheorem*{lemma1.3}{Theorem 1}

\newtheorem*{lemma2}{Lemma 2}

\newtheorem*{lemma3.1b}{Lemma 3.1b (Unfolding)}
\newtheorem*{lemma3.1a}{Lemma 3.1a (Contour Shift)}
\newtheorem*{proposition3.1}{Proposition 3.1}
\newtheorem*{proposition3.2}{Proposition 3.2}

\newtheorem*{lemma3.3a}{Lemma 3.3a}
\newtheorem*{proposition3.3}{Proposition 3.3}

\newtheorem*{note4.2}{Note}

\newtheorem*{lemma2.4}{Lemma A.4}
\newtheorem*{theoremA.5}{Theorem B}

\begin{document}
\title{Sieving for mass equidistribution} 
\author{Roman Holowinsky}
\address{Department of Mathematics, The Ohio State University, 100 Math Tower, 231 West 18th Avenue, Columbus, OH 43210-1174, USA}
\email{romanh@math.ohio-state.edu}
\thanks{This material is based upon work supported by the National Science Foundation and NSERC and written while at the Fields Institute and University of Toronto.}
\begin{abstract}
We approach the holomorphic analogue to the Quantum Unique Ergodicity conjecture through an application of the Large Sieve. We deal with shifted convolution sums as in [Ho], with various simplifications in our analysis due to the knowledge of the Ramanujan-Petersson conjecture in this holomorphic case.
\end{abstract}

\maketitle

\section{Introduction and statement of results}
We study the \textit{shifted convolution sums}
\begin{equation*}
\sum_{n\leqslant x} \lambda_1(n) \overline{\lambda_2(n+\ell)}
\end{equation*}
where $\ell$ is a fixed non-zero integer and the $\lambda_i$ are multiplicative functions.  Ignoring cancellations among the summation terms, we look to obtain non-trivial upper bounds for the sums
\begin{equation*}
\sum_{n \leqslant x} | \lambda_1 ( n ) \lambda_2 ( n + \ell ) |
\end{equation*} 
when the multiplicative functions satisfy $|\lambda_i(n)|\leqslant \tau_m(n)$ for some $m$.  Here $\tau_m(n)$ is the number of ways to represent $n$ as the product of $m$ natural numbers.  This is a continuation of the work in [Ho] where we studied such sums with Hecke eigenvalues of Hecke-Maass cusp forms. The results obtained here will be specifically applied to the holomorphic cusp form analogue since the eigenvalues in this case satisfy the Ramanujan-Petersson bound $|\lambda_f(n)| \leqslant \tau(n) := \tau_2(n)$. 

Our work here complements the recent work of Soundararajan [So] and in [H-S] we combine our results to resolve the mass equidistribution conjecture, the motivation for this work which we now describe.

\subsection{Mass Equidistribution}
Let $\mathbbm{H}$ be the upper half plane with hyperbolic measure $d \mu z : =
y^{- 2} dxdy$.  Set $\Gamma ={\rm SL}\sb 2(\mathbbm{Z})$ and let $X = \Gamma
\backslash \mathbbm{H}$. Denote by $\mathcal{L}( X )$ the Hilbert space of square integrable
automorphic functions with inner product
\begin{equation*} 
\langle f, g \rangle : = \int_X f ( z ) \overline{g ( z )} d \mu z. 
\end{equation*}
Choosing any sequence of holomorphic Hecke eigencuspforms $\{f_k\}$ normalized so that
\begin{equation*}
\int_X y^k |f_k(z)|^2 d\mu z = 1
\end{equation*}
and any fixed $\phi\in \mathcal{L}(X)$, we are interested in the behavior of $\langle \phi F_k,F_k\rangle$ as $k \longrightarrow \infty$.  Here we have simply taken $F_k(z)=y^{k/2} f_k(z)$ where $z=x+iy$.  In particular, we are motivated by the following conjecture.
\begin{Conjecture}
For any fixed $\phi$ smooth and bounded on $X$ we have
\begin{eqnarray}
  \langle \phi F_k, F_k \rangle = \frac{1}{\textnormal{Vol}(X)} \langle \phi, 1 \rangle + o ( 1 ) & \textnormal{as} & k \longrightarrow
  \infty.  \label{1}
\end{eqnarray}
\end{Conjecture}
\noindent Here, $\textnormal{Vol}(X):= \int_X d\mu z = \pi/3$.  Note that this statement is trivial for $\phi$ a constant function.

This equidistribution statement is a holomorphic analogue to a more natural conjecture known as Quantum Unique Ergodicity(QUE) which was originally stated for compact arithmetic surfaces by Rudnick and Sarnak [R-S] and proven in the compact case by Lindenstrauss [Li]. Techniques applied by Lindenstrauss from Ergodic theory for QUE, however, do not seem to carry over to the analogous holomorphic statement of mass equidistribution.  

The statement of QUE, for our choice of non-compact $X$, replaces the varying normalized holomorphic forms in (\ref{1}) with an orthonormal basis of Hecke-Maass cusp forms $\{u_j\}$ and affirms that
\begin{eqnarray*}
  \langle \phi u_j, u_j \rangle = \frac{3}{\pi} \langle \phi, 1 \rangle + o ( 1 ) & \textnormal{as} & t_j \longrightarrow
  \infty.
\end{eqnarray*}
The $u_j$ above are indexed with respect to their Laplace eigenvalues 
\begin{eqnarray*}
  \left( \Delta + (\frac{1}{4}+t_j^2) \right) u_j = 0
\end{eqnarray*}
where
\begin{equation*}
\Delta : = y^2 ( \frac{\partial^2}{\partial x^2} +
   \frac{\partial^2}{\partial y^2} ) 
\end{equation*}
is the hyperbolic Laplacian. The study of the inner products $\langle \phi u_j, u_j \rangle$ as $t_j\longrightarrow\infty$ provides connections between a classically chaotic dynamical system, Hamiltonian flow on a surface of negative curvature, and its quantum model (see \cite{Sa} and \cite{Ze}). 
\subsection{Main results}
We will be analyzing the inner products seen in the statement of mass equidistribution with test functions $\phi$ which have a Fourier series expansion of the form
\begin{equation*} 
\phi ( z ) = \sum_{\ell} a_{\ell} ( y ) e ( \ell x ).
\end{equation*}
For each Fourier coefficient $a_\ell(y)$, we will define the integral
\begin{equation*}
S_\ell(Y) :=  \int^{\infty}_0 g ( Yy ) y^{- 2} \left( \int^{1 / 2}_{- 1
   / 2}\{a_\ell(y) e(\ell x)\} |F_k ( z ) |^2 dx \right) dy
\end{equation*}
from which the shifted convolution sums will arise upon integration in the variable $x$. Here we have chosen $g\in C^\infty_c(\mathbbm{R}^+)$ to be some smooth, compactly supported and positive valued test function and $Y\geqslant 1$ to be some parameter.  The role of $g$ and $Y$ will be to control the size of the shifts $\ell$ which must be considered. These parameters will be introduced in conjunction with an incomplete Eisenstein series
\begin{equation*}
 E_Y(z|g) := \sum_{\gamma \in \Gamma_{\infty} \backslash \Gamma} g (Y \textnormal{Im} (\gamma z) )
\end{equation*} 
which will act as the main ``analytic tool" in manipulating the inner products $\langle \phi F_k , F_k \rangle$.  In fact, our analysis in \S 3 starts with the inner product $\langle E_Y(z|g) \phi F_k, F_k \rangle$.  The extra incomplete Eisenstein series will allow for an ``unfolding" method (see \S 3.1) which will bring us to objects of the form $S_\ell(Y)$ and allows us to work in a Siegel domain with $y\gg 1/Y$ by the support of $g$. Note that in the work of [L-S], the test functions $\phi$ are chosen to be incomplete Eisenstein and Poincare series. With such test functions $\phi$, the unfolding of $\langle \phi F_k , F_k \rangle$ is automatic, but there is less control over the size of the shifts $\ell$.  In \S 3 we prove the following.
\begin{theorem1}
Let $Y \geqslant 1$ and $\varepsilon>0$. Fix a positive valued test function $g\in C^\infty_c(\mathbbm{R}^+)$ and set $c_Y:=\frac{3}{\pi} \langle E(z|g), 1 \rangle Y$. Fix an automorphic form $\phi$ with Fourier expansion
\begin{equation*} 
\phi ( z ) = \sum_{\ell} a_{\ell} ( y ) e ( \ell x ).
\end{equation*}
If $\phi$ is a Hecke-Maass cusp form, then 
\begin{equation*}
\langle \phi F_k, F_k \rangle  = c_Y^{-1} \sum_{0<|\ell|<Y^{1+\varepsilon}} S_\ell(Y) +O\left(\frac{1}{Y^{1/2}}\right).
\end{equation*}
 If $\phi$ is an incomplete Eisenstein series, then
\begin{equation*}
\langle \phi F_k, F_k \rangle = \frac{3}{\pi}\langle \phi,1 \rangle + c_Y^{-1} \sum_{0<|\ell|<Y^{1+\varepsilon}} S_\ell(Y) + O\left(\frac{1+R_k(f)}{Y^{1/2}}\right)
\end{equation*}
with
\begin{equation*}
R_k(f) = \frac{1}{k^{1/2} L(1, \textnormal{sym}^2 f_k)}\int^{+\infty}_{-\infty}\frac{|L(\frac{1}{2}+it, \textnormal{sym}^2 f_k)|}{(|t|+1)^{10}} |dt|.
\end{equation*}
Furthermore, we have the bound
\begin{equation*}
c_Y^{-1} S_\ell(Y) \ll\frac{|a_\ell(Y^{-1})|}{L(1, \textnormal{sym}^2 f_k)}\left\{ \frac{1}{Yk}\sum_{n} |\lambda_f ( n )\lambda_f(n+\ell)| \, g\left(\frac{Y(k-1)}{4 \pi (n+\frac{\ell}{2})}\right) + \frac{(Yk)^{\varepsilon}}{k}\right\}.  
\end{equation*}
\end{theorem1}

We therefore see the interest in the study of the shifted convolution sums which appear in our bound for $c_Y^{-1} S_\ell(Y)$ above.  In \S 4, we prove the following.
\begin{theorem2}
Let $\lambda_1(n)$ and $\lambda_2(n)$ be multiplicative functions satisfying\\ $|\lambda_i(n)|\leqslant\tau_m(n)$ for some $m$.  For any $0<\varepsilon<1$, any $x$ sufficiently large with respect to $\varepsilon$ and any fixed integer $0<|\ell|\leqslant x$ we have
\begin{equation*} 
\sum_{n \leqslant x} | \lambda_1 ( n ) \lambda_2 ( n + \ell ) | \ll  x (\log x)^\varepsilon  M(x) \tau(|\ell|)
\end{equation*}
where
\begin{equation*}
M(x):=\frac{1}{(\log x)^2}\prod_{p \leqslant z} \left(1+\frac{|\lambda_1(p)|}{p}\right)\left(1+\frac{|\lambda_2(p)|}{p}\right)
\end{equation*}
with $z=\exp(\frac{\log x}{\varepsilon \log \log x})$.
\end{theorem2}
\begin{Breteche}
\textnormal{For the sake of completeness, we provide a direct proof of Theorem 2 using the Large Sieve.  However, an alternate proof may be derived from the works of Nair [N] and Nair-Tenenbaum [N-T].}
\end{Breteche}
Using Theorem 2 with $\lambda_1=\lambda_2=\lambda_f$ for each Hecke eigencuspform $f_k$ and bounding $|a_\ell(Y^{-1})|$ in Theorem 1 by an application of Lemma 2 below (with $A=0$), leaves us with evaluating a simple sum over the shifts $0<|\ell|<Y^{1+\varepsilon}$.  In the cusp form case we are left with
\begin{equation*}
 \frac{(Y \log k)^{\varepsilon}}{Y^{1/2}(\log k)^2 L(1, \textnormal{sym}^2 f_k)} \prod_{p\leqslant k} \left(1+\frac{2 |\lambda_f(p)|}{p}\right) \sum_{0<|\ell|<Y^{1+\varepsilon}} |\rho(\ell)| \tau(|\ell|).
\end{equation*} 
In the incomplete Eisenstein series case, $|\rho(\ell)|$ is replaced by $\tau(|\ell|)$.  Cauchy's inequality and the bound \eqref{coeffsquaredavg} below (used only in the cusp form case) gives the following Corollary upon choosing $Y$ optimally.

\begin{corollary3}
Fix an automorphic form $\phi$. Define
\begin{equation*}
M_k(f) := \frac{1}{(\log k)^2 L(1, \textnormal{sym}^2 f_k)} \prod_{p\leqslant k} \left(1+\frac{2 |\lambda_f(p)|}{p}\right). 
\end{equation*}
If $\phi$ is a Hecke-Maass cusp form, then 
\begin{equation*}
\langle \phi F_k, F_k \rangle  \ll (\log k)^{\varepsilon} M_k(f)^{1/2}
\end{equation*}
for any $\varepsilon>0$.  If $\phi$ is an incomplete Eisenstein series then
\begin{equation*}
\langle \phi F_k, F_k \rangle = \frac{3}{\pi}\langle \phi,1 \rangle + O\left((\log k)^{\varepsilon} M_k(f)^{1/2}(1 + R_k(f))\right)
\end{equation*}
for any $\varepsilon>0$ with
\begin{equation*}
R_k(f) \ll \frac{1}{k^{1/2} L(1, \textnormal{sym}^2 f_k)}\int^{+\infty}_{-\infty}\frac{|L(\frac{1}{2}+it, \textnormal{sym}^2 f_k)|}{(|t|+1)^{10}} |dt|.
\end{equation*}
\end{corollary3}
\begin{note1}
\textnormal{In Theorem 1 we ask that $Y\geqslant 1$.  The bounds in Corollary 3 come from choosing $Y = M_k(f)^{-1}$.   We trivially have that $M_k(f)^{\varepsilon} \ll (\log k)^{\varepsilon}$. If $M_k(f)>1$ for some form $f_k$, then simply choosing $Y=1$ in Theorem 1 (i.e. ignoring the shifted sums) produces the better result.  In either case, the bounds in Corollary 3 still hold true.}  
\end{note1} 

Applying the Ramanujan-Petersson bound trivially in Corollary 3, one fails to properly control the growth rate of $M_k(f)$. The usefulness of Theorem 2 might therefore initially be unapparent. However, we still expect the main term $M_k(f)$ to be relatively small due to the conjectured Sato-Tate distribution. More specifically, we expect that
\begin{equation*}
\exp\left(2\sum_{p\leqslant k} \frac{|\lambda_f(p)|-1}{p} \right) = o(L(1, \textnormal{sym}^2 f_k))
\end{equation*}
as $k\longrightarrow \infty$.  To help illustrate this point, we turn to an idea seen in [E-M-S]. Simply by the Ramanujan-Petersson bound we have
\begin{equation*}
2|\lambda_f(p)|-2 \leqslant (\lambda_f^2(p)-1) -\frac{1}{9}  (\lambda_f^2(p)-1)^2
\end{equation*}
for each prime $p$. The Hecke relations for the coefficients of the symmetric square then give that
\begin{equation*}
\sum_{p\leqslant k} \frac{2|\lambda_f(p)|-2}{p} \leqslant \sum_{p\leqslant k} \frac{\lambda_f(p^2)}{p} -\frac{1}{9} \sum_{p\leqslant k} \frac{\lambda_f^2(p^2)}{p}.
\end{equation*}
We expect the contribution from the sum over terms of the form $\lambda(p^m) p^{-1}$ to be roughly of size $\log L(1, \textnormal{sym}^m f_k)$ and with $\lambda_f^2(p^2) = 1+\lambda_f(p^2)+\lambda_f(p^4)$ one therefore hopes to obtain a bound of the form
\begin{equation*}
M_k(f) \ll \{ (\log k) L(1, \textnormal{sym}^2 f_k)  L(1, \textnormal{sym}^4 f_k)\}^{-1/9} \ll_f (\log k)^{-\delta}
\end{equation*}
for some $\delta>0$. Recall that the holomorphicity and non-vanishing of the symmetric $m$-th power $L$-functions in the half-plane $\textnormal{Re}(s)\geqslant 1$ is known for $m\leqslant 8$ by the works of [K-Sh] and [Ki].  Although it seems difficult, by means of current technology, to show that the above bound for $M_k(f)$ is uniform for all $f$, one may at least show that the above bound holds uniformly for some $\delta>0$ for all but at most $O(k^\varepsilon)$ forms of weight $k$ by a zero-free region argument.

In his recent work [So], Soundararajan attacks the problem of mass equidistribution through the analysis of $L$-functions at their central values.  Soundararajan's method also produces mass equidistribution with a small number of possible exceptions. In [H-S], we combine our two different arguments to remove the possibility of any exceptions thus proving mass equidistribution with an effective rate of convergence. 

\begin{acknowledgements}
\textnormal{The author thanks all of those who helped contribute to the development of this work through useful comments, suggestions and guidance.  In particular, thanks go to Prof. H. Iwaniec, Prof. P. Sarnak, and Prof. K. Soundararajan. The author also thanks Prof. R. de la Bret\`{e}che for pointing out an alternate proof of Theorem 2 using the works of Nair [N] and Nair-Tenenbaum [N-T].}
\end{acknowledgements}

\section{Fourier coefficients of the forms $\phi$ and $f_k$}
For $\phi$ a fixed Maass cusp form or incomplete Eisenstein series, we can express $\phi$
as a Fourier series expansion of the type
\begin{equation} \phi ( z ) = a_0 ( y ) + \sum_{\ell \neq 0} a_{\ell} ( y ) e ( \ell x ), \label{2.3} \end{equation}
with $a_0 ( y )= \langle\phi,1\rangle = 0$ in the cusp form case.  If $\phi ( z )$ is a fixed Maass cusp
form with $\Delta$ eigenvalue $1/4+r^{2}$, then we have the expansion
\begin{equation*} \phi ( z ) = \sqrt{y}  \sum_{\ell \neq 0} \rho ( \ell ) K_{ir} ( 2 \pi | \ell
   |y ) e ( \ell x ) \end{equation*}
where the $\rho ( \ell )$ are complex numbers. We know by ([Iw], Thrm 3.2) that the Fourier coefficients of a cusp form satisfy
\begin{equation}
\sum_{|\ell|\leqslant x} |\rho(\ell)|^2 \ll (r + x) e^{\pi r}. \label{coeffsquaredavg}
\end{equation}

For $\phi ( z ) = E ( z| \psi )$
the incomplete Eisenstein series
\begin{equation*} 
E ( z| \psi ) : = \sum_{\gamma \in \Gamma_{\infty} \backslash \Gamma} \psi (
   \textnormal{Im } \gamma z ), 
\end{equation*}
where $\psi ( y )$ is a smooth function, compactly supported on
$\mathbbm{R}^+$, the coefficients $a_{\ell} ( y )$ can be determined in terms
of the coefficients of the Eisenstein series
\begin{equation} 
E ( z, s ) : = \sum_{\gamma \in \Gamma_{\infty} \backslash \Gamma} ( \textnormal{Im }
   \gamma z )^s \label{Eisenstein} . 
\end{equation}
The latter has the Fourier expansion
\begin{equation} E ( z, s ) = y^s + \varphi ( s ) y^{1 - s} + \sqrt{y}  \sum_{\ell \neq 0}
   \varphi_{| \ell |} ( s ) K_{s - \frac{1}{2}} ( 2 \pi | \ell |y ) e ( \ell x
   ) \label{2.4} \end{equation}
where
\begin{eqnarray*}
  \varphi ( s ) = & \sqrt{\pi}  \frac{\Gamma ( s - \frac{1}{2} ) \zeta ( 2 s -
  1 )}{\Gamma ( s ) \zeta ( 2 s )} = & \frac{\theta ( 1 - s )}{\theta ( s
  )},\\
  \theta ( s ) = & \pi^{- s} \Gamma ( s ) \zeta ( 2 s ), & \\
  \varphi_{\ell} ( s ) = & \frac{2}{\theta ( s )}  \sum_{ab = \ell} \left(
  \frac{a}{b} \right)^{s - \frac{1}{2}}, & \textnormal{if } \ell \geqslant 1.
\end{eqnarray*}
Note that
\begin{equation}
  \textnormal{Res}_{s = 1} E ( z, s ) = \textnormal{Res}_{s = 1} \varphi ( s ) = \frac{3}{\pi} . \label{2.5}
\end{equation}
We have
\begin{equation*} E ( z| \psi ) = \frac{1}{2 \pi i} \int_{( 2 )} \Psi ( - s ) E ( z, s ) ds 
 \end{equation*}
where $\Psi ( s )$ is the Mellin transform of $\psi ( y )$.  This is an entire
function with rapid decay in vertical strips, specifically
\begin{equation*} \Psi ( s ) \ll ( |s| + 1 )^{- A} \end{equation*}
for any $A \geqslant 0$, $- 2 \leqslant \textnormal{Re} ( s ) \leqslant 2$, with the
implied constant depending only on $\psi$ and $A$.

From these formulas, we find the coefficients of (\ref{2.3}) in the incomplete
Eisenstein series case
\begin{equation*} a_0 ( y ) = \frac{1}{2 \pi i} \int_{( 2 )} \Psi ( - s ) ( y^s + \varphi ( s
   ) y^{1 - s} ) ds = \psi ( y ) + O ( y^{- 1} ) \end{equation*}
and for $\ell \neq 0$ we move the integration to the line $\textnormal{Re}(s) = 1/2$ to get
\begin{equation} a_{\ell} ( y ) = \left( \frac{y}{\pi} \right)^{\frac{1}{2}} \int^{+
   \infty}_{- \infty} \frac{\pi^{it} \Psi ( - \frac{1}{2} - it )}{\Gamma (
   \frac{1}{2} + it ) \zeta ( 1 + 2 it )}  \left( \sum_{ab = |\ell|} \left(
   \frac{a}{b} \right)^{it} \right) K_{it} ( 2 \pi | \ell |y ) dt. \label{2.6}
\end{equation}
Doing the same for $a_0 ( y )$ we get
\begin{equation*} a_0 ( y ) = \frac{3}{\pi} \Psi ( - 1 ) + O ( \sqrt{y} ). \end{equation*}
On the other hand, by unfolding the incomplete Eisenstein series $\phi$, we
derive that
\begin{equation} \langle \phi, 1 \rangle = \int^{1 / 2}_{- 1 / 2} \int_0^{\infty} \psi ( y ) d \mu z = \Psi
   ( - 1 ). \label{testconnection}
\end{equation}
Therefore, we have
\begin{equation*} 
a_0 ( y ) = \frac{3}{\pi} \langle \phi, 1 \rangle + O ( \sqrt{y} )  
\end{equation*}
which conveniently contains the expected main term $3/\pi \langle\phi,1\rangle$ in the statement of mass equidistribution.  For $\ell \neq 0$ we apply to (\ref{2.6}) the bound
\begin{equation*} 
K_{it} ( w ) \ll | \Gamma ( \frac{1}{2} + it ) | \left( \frac{1 + |t|}{w} \right)^A \left(1+\frac{1+|t|}{w}\right)^{\varepsilon} 
\end{equation*}
known for real $w$ and imaginary order $it$ with any $\varepsilon>0$ and any integer $A\geqslant 0$ by repeated integration by parts of the integral representation of the $K$-Bessel function
\begin{equation*} 
K_{it} ( w ) = \pi^{- 1 / 2} \Gamma \left(\frac{1}{2} + it \right)
   \left( \frac{w}{2} \right)^{- it}  \int_0^{+ \infty} ( v^2 + 1 )^{- it -
   1 / 2} \cos ( vw ) dv ,
\end{equation*}
in order to obtain
\begin{equation*} 
a_{\ell} ( y ) \ll \tau ( | \ell | ) \sqrt{y} \left(\frac{1}{|\ell|y} \right)^A \left(1+\frac{1}{|\ell|y}\right)^{\varepsilon}. 
\end{equation*}
Here $\tau ( \ell )$ is the divisor function. In the case of the cusp form, when $a_\ell(y) = \rho(\ell) K_{ir}(2 \pi|\ell| y)$, we
get a similar bound. We state the results in the following Lemma.
\begin{lemma2}
  Let $\phi \in \mathcal{A} ( X )$ be an automorphic function on $X =
 {\rm SL}\sb 2(\mathbbm{Z}) \backslash \mathbbm{H}$ with Fourier series expansion
  \begin{equation*} \phi ( z ) = a_0 ( y ) + \sum_{\ell \neq 0} a_{\ell} ( y ) e ( \ell x ) . \end{equation*}
 If $\phi$ is a Maass cusp form with $\Delta$ eigenvalue $1/4+r^{2}$, then $a_0 ( y ) = 0$ and for $\ell \neq 0$ we have 
  \begin{equation*}
    a_{\ell} ( y ) \ll | \rho ( \ell ) |\sqrt{y} \left(\frac{1+|r|}{|\ell|y} \right)^A \left(1+\frac{1+|r|}{|\ell|y}\right)^{\varepsilon}
  \end{equation*}
for any $\varepsilon>0$ and any integer $A\geqslant 0$.  If $\phi$ is an incomplete Eisenstein series, then 
\begin{equation*} 
a_0 ( y )  =  \frac{3}{\pi} \langle \phi, 1 \rangle + O ( \sqrt{y} )
\end{equation*} 
and for $\ell \neq 0$ we have
  \begin{equation*}
   a_{\ell} ( y ) \ll  \tau ( | \ell | ) \sqrt{y} \left(\frac{1}{|\ell|y} \right)^A \left(1+\frac{1}{|\ell|y}\right)^{\varepsilon}
  \end{equation*}
for any $\varepsilon>0$ and any integer $A\geqslant 0$.
\end{lemma2}
For $f_k$ a holomorphic Hecke eigencuspform of weight $k$ with Petersson norm
\begin{equation*}
\int_{X}y^k|f_k(z)|^2d\mu z = 1,
\end{equation*}
we can express $f_k$ as a Fourier series expansion of the form
\begin{equation*} 
f_k ( z ) = \sum_{n \geqslant 1} a_{f} ( n ) e ( n z ). 
\end{equation*}
The coefficients $a_f(n)$ are proportional to the Hecke eigenvalues $\lambda_f(n)$
\begin{equation}
a_f(n)=\lambda_f(n)a_f(1) n^{(k-1)/2} \label{fourier2hecke}
\end{equation}
with the first Fourier coefficient satisfying
\begin{equation}
|a_f(1)|^2=\frac{(4\pi)^{k-1}}{\Gamma(k-1)}\frac{2 \pi^2}{(k-1)L(1,\textnormal{sym}^2 f_k)}.\label{firstfouriersquared}
\end{equation}
The Ramanujan-Petersson conjecture, $|\lambda_f(p)|\leqslant 2$ for prime $p$, is known for holomorphic forms $f_k$ by the work of Deligne and it is also known ([H-L], [G-H-L]) that
\begin{equation*}
(\log k)^{-1} \ll L(1, \textnormal{sym}^2 f_k)\ll (\log k)^3.
\end{equation*} 
\section{Proof of Theorem 1}
We demonstrate how one may relate the inner products 
\begin{equation}
 < \phi F_k, F_k >=\int_{X}\phi(z) |F_k(z)|^2 d \mu z \label{maininner}
\end{equation}
in the question of mass equidistribution to the study of a controlled number of shifted convolution sums.  As mentioned in the introduction, we will introduce an additional incomplete Eisenstein series along with a parameter $Y$ which will act as our main analytic tools.
\subsection{Construction of our main object $I_\phi(Y)$}
Let $f_k$ be a normalized holomorphic Hecke eigencuspform of weight $k$ and $\phi$ an automorphic function which is smooth and bounded on $\mathbbm{H}$. Let $F_k(z):=y^{k/2}f_k(z)$.  Fix a function $g ( y ) \in
C^{\infty}_c (\mathbbm{R}^+ )$, smooth and compactly supported on
$\mathbbm{R}^+$, and let
\begin{equation*} G ( s ) := \int^{+ \infty}_0 g ( y ) y^{s - 1} dy \end{equation*}
be its Mellin transform.  Therefore, $G ( s )$ is entire and
\begin{equation*} G ( s ) \ll ( |s| + 1 )^{- A} \end{equation*}
for any $A \geqslant 0$, uniformly in vertical strips, where the implied
constant depends on $g$ and $A$. Let $Y \geqslant 1$ and consider the integral
\begin{equation}
 I_{\phi} ( Y ) := \frac{1}{2 \pi i} \int_{( \sigma )} G ( - s ) Y^s
   \int_{X} E ( z, s ) \phi ( z ) |F_k ( z ) |^2 d \mu zds \label{3.2} \end{equation}
with $\sigma > 1$.   The integral $I_\phi(Y)$ will provide the connection between our inner products (\ref{maininner}) and shifted convolution sums.  We compute $I_\phi(Y)$ asymptotically in two ways.
\begin{lemma3.1a} 
For $\phi$ a fixed Hecke-Maass cusp form or incomplete Eisenstein series we have
\begin{equation*}
  I_\phi ( Y ) = c_Y \langle \phi F_k, F_k \rangle  + O ( Y^{1/2} ) 
\end{equation*}
where 
\begin{equation}
c_Y:=\frac{3}{\pi} \langle E(z|g), 1 \rangle Y.  \label{C}
\end{equation}
\begin{proof}
Starting with equation (\ref{3.2}) and moving the contour of
integration to the line $\textnormal{Re}(s) = 1/2$, we write
\begin{equation*}
  I_\phi ( Y ) = c_Y \langle \phi F_k, F_k \rangle + R_\phi ( Y )
\end{equation*}
with $c_Y$ as in \eqref{C} coming from the pole of the Eisenstein series at $s = 1$ (see \eqref{testconnection}) and $R_\phi ( Y )$ the remaining term
\begin{equation*}
  R_\phi ( Y ): =  \int_{X} p(z) \phi ( z ) |F_k ( z ) |^2 d
  \mu z 
\end{equation*}
where
\begin{equation*}
p(z):=\frac{1}{2 \pi i} \int_{( 1 / 2 )} G ( - s ) Y^s E ( z, s ) ds.
\end{equation*}
From the Fourier series expansion (\ref{2.4}) for the Eisenstein series $E ( z, s )$
and Lemma 2 (with $A=2$) we have $E (z,s) \ll  \sqrt{y} + |s|^2 y^{- 3 / 2}(1+|s|/y)^{\varepsilon}$ on the line $\textnormal{Re} ( s ) = 1 / 2$. Hence $p ( z ) \ll \sqrt{yY}$ if $y
\geqslant 1 / 2$.  Assuming that $\sqrt{y} |\phi ( z ) |$ is bounded on
$\mathbbm{H}$, we conclude that $R_\phi ( Y ) \ll_{\phi,g} \sqrt{Y}$.  Such is the case for $\phi$ a cusp form or incomplete Eisenstein series.
\end{proof}
\end{lemma3.1a}
\begin{lemma3.1b}
For $\phi$ a Hecke-Maass cusp form or incomplete Eisenstein series we have
\begin{equation} 
I_{\phi} ( Y )=\int^{\infty}_0 g ( Yy ) y^{- 2} \left( \int^{1 / 2}_{- 1 / 2} \phi ( z
  ) |F_k ( z ) |^2 dx \right) dy . \label{lemma2.1} 
\end{equation}
\begin{proof}
Unfolding the inner integral in (\ref{3.2}) by the definition of the Eisenstein series \eqref{Eisenstein} and then integrating in $s$ we get the result. 
\end{proof}
\end{lemma3.1b}
\begin{remark0}
\textnormal{Since $g ( Yy )$ is supported on $y \asymp 1 /
Y$, $|x|\leqslant 1/2$ and $\phi ( z )$ is bounded, we have by Lemma 2.10 in [Iw], which states that there are roughly $Y$ copies of the standard fundamental domain in that region, that $I_{\phi} ( Y ) \ll_{\phi, g} Y$. Note that the same bound holds for $\phi \equiv 1$ the constant function which we use in \S 3.2.}
\end{remark0} 
\noindent Combining Lemma 3.1a with Lemma 3.1b we get the following.
\begin{proposition3.1}
For $\phi$ a Hecke-Maass cusp form or incomplete Eisenstein series we have
\begin{equation*}
 c_Y \langle \phi F_k, F_k \rangle  + O ( Y^{1/2} ) =  \int^{\infty}_0 g ( Yy ) y^{- 2} \left( \int^{1 / 2}_{- 1 / 2} \phi ( z
  ) |F_k ( z ) |^2 dx \right) dy 
\end{equation*}
where $c_Y= \frac{3}{\pi} \langle E(z|g), 1 \rangle Y$.
\end{proposition3.1}
\noindent The remainder of \S 3 will be devoted to the analysis of the integral in \eqref{lemma2.1} to provide a useful relation between $\langle \phi F_k, F_k \rangle$ and shifted convolution sums via Proposition 3.1.
\subsection{Truncating the test function $\phi(z)$}
We are left with extracting shifted convolution sums from the integral 
\begin{equation*}
I_{\phi} ( Y )=\int^{\infty}_0 g ( Yy ) y^{- 2} \left( \int^{1 / 2}_{- 1 / 2} \phi ( z
  ) |F_k ( z ) |^2 dx \right) dy.
\end{equation*}
We first choose to truncate our fixed form $\phi$ in order to restrict the size of the shifts $\ell$ which must be considered.  Recall that the fixed form $\phi ( z )$ has Fourier expansion
\begin{equation} 
\phi ( z ) = \sum_{\ell} a_{\ell} ( y ) e ( \ell x ) \label{4.111111111} 
\end{equation}
with coefficient bounds as seen in Lemma 2.  We ignore the dependence on the spectral parameter for bounds involving $\phi$ since it is fixed in this case. 

If $\phi$ is an incomplete Eisenstein series, then we find that the contribution to $I_\phi(Y)$ from the tail of (\ref{4.111111111}) with $| \ell | \geqslant Y^{1+\varepsilon}$ for any $\varepsilon>0$ is bounded by
\begin{equation*} 
\left(\int^{\infty}_0 \int^{1 / 2}_{- 1 / 2} g(Yy) |F_k ( z ) |^2 \frac{dx dy}{y^2} \right) \frac{Y^{A+\varepsilon}}{\sqrt{Y}}\sum_{| \ell | \geqslant Y^{1+\varepsilon}} \frac{\tau ( \ell )}{|\ell|^A} \ll Y^{3/2+\varepsilon(1-A)} 
\end{equation*}
by the support of $g$ and Lemma 2.  See the remark after Lemma 3.1b.  Hence the contribution of these terms to $I_{\phi} ( Y )$ is $\ll Y^{1/2}$ after choosing $A$ sufficiently large with respect to $\varepsilon$. A similar argument works when $\phi$ is a fixed cusp form.   
\begin{proposition3.2}
For $\phi$ a Hecke-Maass cusp form or incomplete Eisenstein series we have
\begin{equation*} 
 c_Y \langle \phi F_k, F_k \rangle =  \int^{\infty}_0 g ( Yy ) y^{- 2} \left( \int^{1 / 2}_{- 1 / 2} \phi^{\ast} ( z
  ) |F_k ( z ) |^2 dx \right) dy + O(Y^{1/2}) 
\end{equation*}
where $c_Y= \frac{3}{\pi} \langle E(z|g), 1 \rangle Y$ and
\begin{equation*} 
\phi^{\ast} ( z ) := \sum_{| \ell | < Y^{1+\varepsilon}} a_{\ell} ( y ) e ( \ell x ). 
\end{equation*}
\end{proposition3.2}
\subsection{Extracting shifted convolution sums}
We now consider each Fourier coefficient $a_\ell(y)$ of $\phi^{\ast}(z)$ separately and write 
\begin{equation} 
c_Y \langle \phi F_k, F_k \rangle =S_0(Y) + \sum_{0<|\ell|< Y^{1+\varepsilon}} S_{\ell} ( Y ) + O(Y^{1/2})\label{4.2}
\end{equation}
where for any integer $\ell$ we define
\begin{equation}
S_\ell(Y) :=  \int^{\infty}_0 g ( Yy ) y^{- 2} \left( \int^{1 / 2}_{- 1
   / 2}\{a_\ell(y) e(\ell x)\} |F_k ( z ) |^2 dx \right) dy. \label{Sell}
\end{equation}
The aim of this section is to analyze the objects $S_\ell(Y)$ so that when we divide through equation \eqref{4.2} by $c_Y$ we will have the equations and bounds seen in Theorem 1.  We start by noting that $S_0(Y)\equiv 0$ for $\phi$ a cusp form and by Lemma 2 we have
\begin{equation}
S_0(Y) =  \{\frac{3}{\pi}\langle\phi,1\rangle + O(Y^{-1/2})\} \int^{\infty}_0 g ( Yy ) y^{- 2} \left( \int^{1 / 2}_{- 1 / 2} |F_k ( z ) |^2 dx \right) dy \label{S0}
\end{equation}
for $\phi$ an incomplete Eisenstein series. We treat $S_0(Y)$  in \eqref{S0} and $S_\ell(Y)$ in \eqref{Sell} for $\ell \neq 0$ simultaneously.  

Squaring out $|F_k(z)|^2$ and integrating in $x$ gives
\begin{equation*}
S_0(Y) =\{\frac{3}{\pi}\langle\phi,1\rangle + O(Y^{-1/2})\} \sum_{n \geqslant 1} 
|a_f ( n )|^2 \left( \int_0^{\infty} g(Yy) y^{k-2} e^{-4\pi n y}dy\right)
\end{equation*}
and
\begin{equation*}
  S_\ell( Y ) =  \sum_{n \geqslant 1} 
\overline{a_f ( n )} a_f(n+\ell) \left( \int_0^{\infty} g(Yy) a_\ell(y) y^{k-2} e^{-2\pi (2n+\ell) y}dy\right)
\end{equation*}
for $\ell \neq 0$ which satisfies
\begin{equation*}
S_\ell(Y) \ll |a_\ell(Y^{-1})| \sum_{n\geqslant 1} |a_f(n) a_f(n+\ell)| \left( \int_0^{\infty} g(Yy) y^{k-2} e^{-2\pi (2n+\ell) y}dy\right)
\end{equation*}
by choice of test function $g$ taking non-negative values.  Appealing to the Mellin transform of $g(Yy)$ and applying the normalizations (\ref{fourier2hecke}) and (\ref{firstfouriersquared}) for $a_f(n)$ and $|a_f(1)|^2$, we integrate in $y$ and define for all integers $n$ and $\ell$ with $n\geqslant 1$ and $n+\ell \geqslant 1$ 
\begin{equation*}
W_{n,\ell}(Y):=\left(\frac{\sqrt{n(n+\ell)}}{n+\frac{\ell}{2}}\right)^{k-1}\frac{1}{2 \pi i} \int_{(\sigma)} G(-s) \left(\frac{Y}{4 \pi (n+\frac{\ell}{2})}\right)^s \frac{\Gamma(s+k-1)}{\Gamma(k-1)} ds 
\end{equation*}
for any $\sigma>1$.  Note that $\sqrt{n(n+\ell)}\leqslant n+\ell/2$ with equality when $\ell=0$.
\begin{lemma3.3a}
Let $Y\geqslant 1$. When $\ell=0$ we have $S_0(Y)\equiv 0$ for $\phi$ a cusp form and
\begin{equation*}
S_0(Y) = \{\frac{3}{\pi}\langle\phi,1\rangle + O(Y^{-1/2})\}\frac{2 \pi^2}{(k-1) L(1,\textnormal{sym}^2 f_k)} \sum_n |\lambda_f(n)|^2 W_{n,0}(Y)
\end{equation*}
for $\phi$ an incomplete Eisenstein series.  When $\ell\neq 0$ we have
\begin{equation*}
S_\ell(Y) \ll  \frac{|a_\ell(Y^{-1})|}{k L( 1, \textnormal{sym}^2 f_k)} \sum_{n} |\lambda_f(n)\lambda_f(n+\ell)|W_{n,\ell}(Y)
\end{equation*}
for both such $\phi$.
\end{lemma3.3a}
In the case of the $\ell=0$ shift, we make the extra observation that
\begin{equation*}
\sum_{n} |\lambda_f(n)|^2 W_{n,0}(Y) = \frac{1}{2 \pi i} \int_{(\sigma)} G(-s) L ( s, f_k \otimes f_k ) \left(\frac{Y}{4 \pi} \right)^{s} \frac{\Gamma(s+k-1)}{\Gamma(k-1)} ds
\end{equation*}
where
\begin{equation*} 
L ( s, f_k \otimes f_k ) := \sum_{n \geqslant 1} \lambda^2_f ( n ) n^{-s} 
\end{equation*}
is the Rankin-Selberg convolution. Moving the contour of integration to the $\sigma =1/2$ line we pick up the pole at $s=1$ giving
\begin{equation*}
\sum_{n} |\lambda_f(n)|^2 W_{n,0}(Y) = \frac{G(-1)Y (k-1) L(1, \textnormal{sym}^2 f_k)}{4 \pi \zeta(2)}+E_{1/2}(Y) 
\end{equation*}
where
\begin{equation*}
E_{1/2}(Y):=\frac{1}{2 \pi i} \int_{(1/2)} G(-s) L ( s, f_k \otimes f_k ) \left(\frac{Y}{4 \pi}\right)^s \frac{\Gamma(s+k-1)}{\Gamma(k-1)} ds.
\end{equation*}

We now make use of a Lemma seen in the work of Luo and Sarnak.  By ([L-S], (2.3)) we have
\begin{equation}
\frac{\Gamma(s+k-1)}{\Gamma(k-1)} = (k-1)^s \left(1+O_{a,b}((|s|+1)^2 k^{-1})\right) \label{gammaratio}
\end{equation}
which holds by Stirling's formula for any vertical strip $0<a\leqslant\textnormal{Re}(s)\leqslant b$. Applying this in the case of the zero shift for $E_{1/2}(Y)$ tells us that
\begin{equation*}
E_{1/2}(Y) \ll (Yk)^{1/2}\int^{+\infty}_{-\infty}\frac{|L(\frac{1}{2}+it, \textnormal{sym}^2 f_k)|}{(|t|+1)^{A}} |dt|. 
\end{equation*} 
for any $A>0$.  For the non-zero shifts $\ell$, \eqref{gammaratio} and integration in $s$ with $\sigma=1+\varepsilon$ for some $\varepsilon>0$ gives
\begin{equation}
W_{n,\ell}(Y)  = \left(\frac{\sqrt{n(n+\ell)}}{n+\frac{\ell}{2}}\right)^{k-1} g\left(\frac{Y(k-1)}{4 \pi (n+\frac{\ell}{2})}\right) + O\left(k^{\varepsilon}\left(\frac{Y}{n+\frac{\ell}{2}}\right)^{1+\varepsilon}\right).\label{Wnell}
\end{equation}
Combining the above, we get the following by applying Lemma 3.3a along with a trivial application of Cauchy's inequality and Deligne's bound to control the contribution from the error term in \eqref{Wnell} above.
 \begin{proposition3.3}
Let $Y\geqslant 1$. When $\ell=0$ we have $S_0(Y)\equiv 0$ for $\phi$ a cusp form and
\begin{eqnarray}
\nonumber c_Y ^{-1} S_0(Y) & = & \{\frac{3}{\pi}\langle\phi,1\rangle + O(Y^{-1/2})\}\\
\nonumber & & \times \left\{1+O\left((Yk)^{-1/2}\int^{+\infty}_{-\infty}\frac{|L(\frac{1}{2}+it, \textnormal{sym}^2 f_k)|}{(|t|+1)^{10}L(1, \textnormal{sym}^2 f_k)} |dt|\right)\right\} 
\end{eqnarray}
for $\phi$ an incomplete Eisenstein series. For $\ell\neq 0$ we have for both such $\phi$
\begin{equation*}
c_Y^{-1} S_\ell(Y) \ll \frac{|a_\ell(Y^{-1})|}{L(1, \textnormal{sym}^2 f_k)}\left\{ \frac{1}{Yk}\sum_{n} |\lambda_f ( n )\lambda_f(n+\ell)| \, g\left(\frac{Y(k-1)}{4 \pi (n+\frac{\ell}{2})}\right) + \frac{(Yk)^{\varepsilon}}{k}\right\} 
\end{equation*}
for any $\varepsilon>0$.
\end{proposition3.3}
\noindent With this and Proposition 3.2, we conclude the proof of Theorem 1.
\section{Proof of Theorem 2}
We first rearrange and partition our shifted sums
\begin{equation}
S_{\ell} ( x ) := \sum_{n \leqslant x} | \lambda_1 ( n ) \lambda_2 ( n + \ell ) |,
\label{SCS}
\end{equation}
into objects which may either be treated by elementary methods or by a Large Sieve.  Recall that we are working with  general multiplicative functions satisfying $|\lambda_i(n)|\leqslant\tau_m(n)$ for some $m$.  For future reference, we assume that $0<|\ell|\leqslant x$ and given $0<\varepsilon<1$ we will be working throughout \S 4 with the choice of variables $2\leqslant z\leqslant y \leqslant x$ satisfying
\begin{eqnarray}
\nonumber z &:= & x^{1/s},\\
\nonumber s &:= & \varepsilon \log \log x,\\
\nonumber y &:= & x^{\varepsilon},\\
 x & \geqslant & \exp(\exp(\exp\{(4+m^4)(2\varepsilon)^{-1}\})).\label{choices}
\end{eqnarray}  

\subsection{Factorization of $n$ and $n+\ell$ and partitioning of $S_\ell(x)$}
Factoring $n$ and $n + \ell$ uniquely as
\begin{eqnarray}
  n = ab & \textnormal{and} & n + \ell = a_\ell b_\ell, \label{9.2.1}
\end{eqnarray}
such that for every prime $p$ dividing $n(n+\ell)$,
\begin{eqnarray}
 p| a a_\ell  \Rightarrow p \leqslant z & \textnormal{and} & p|b b_\ell \Rightarrow p>z, \label{9.2.2}
\end{eqnarray}
we partition the sum $S_{\ell} ( x )$ into parts depending on the size of $a$ and $a_\ell$.  We shall refer to the factorizations in \eqref{9.2.1}-\eqref{9.2.2} as $\mathcal{F}_z$.  The partitioning of the sum $S_\ell(x)$ is both convenient and necessary for the Large Sieve application.  We denote by $\mathcal{S}^y ( x )$ the parts where $a$ or $a_\ell$ are greater than $y$, 
\begin{equation}
  \mathcal{S}^y ( x ) :=  \sum_{\substack{
    n \leqslant x\\
    \mathcal{F}_z \textnormal{ and } a > y
  }} | \lambda_1 ( n ) \lambda_2 ( n + \ell ) | + \sum_{\substack{
    n \leqslant x\\
   \mathcal{F}_z \textnormal{ and } a_\ell > y
  }} | \lambda_1 ( n ) \lambda_2 ( n + \ell ) |, \label{SA}
\end{equation}
and the part where both $a$ and $a_\ell$ are less than or equal to $y$ we denote by $\mathcal{S}_y ( x )$ 
\begin{equation}
\mathcal{S}_{y} ( x ) := \sum_{\substack{
    n \leqslant x\\
   \mathcal{F}_z \textnormal{ and } a, a_\ell \leqslant y
  }} | \lambda_1 ( n ) \lambda_2 ( n+\ell )
  | \label{Sy}
\end{equation}
so that
\begin{equation}
S_{\ell}(x)\leqslant \mathcal{S}^y(x)+\mathcal{S}_{y} (x).
\label{2parts}
\end{equation}
In equations \eqref{SA}-\eqref{Sy}, $\mathcal{F}_z$ denotes that we are considering $n$ and $n+\ell$ in terms of their factorizations given by \eqref{9.2.1}-\eqref{9.2.2}.

\subsection{Treating $\mathcal{S}^y(x)$ by elementary methods}
We first handle the part with $a$ or $a_\ell$ large. The technique need only be elementary, because the number of such $a$ and $a_\ell$ with small prime factors should be relatively small.

By H\"{o}lder's inequality and $|\lambda_i(n)|\leqslant \tau_m(n)$ we get
\begin{equation*}
  \mathcal{S}^y( x )  \ll(\sum_{\substack{
    y< a \leqslant x\\
    p|a \Rightarrow p\leqslant z
  }} \sum_{\substack{ b \leqslant x/a \\
                       p|b \Rightarrow p>z}} 1)^{1 / 2} ( \sum_{n \leqslant x}  \tau_m^4 ( n ) )^{1 /
  2}
 \end{equation*}
for any $|\ell|\leqslant x$. We know (see for example [I-K], (1.80)) that 
\begin{equation*}
\sum_{n\leqslant x} \tau_m^4(n) \ll x(\log x)^{m^4 - 1}.
\end{equation*}
By a classical Rankin's method argument ([M-V], Thrm 7.6), the number of integers up to $x$, containing only small prime factors, can be bounded by $\ll x (\log x)^{-A}$ for any $A>0$.  Partial summation and our choice of $2\leqslant z \leqslant y\leqslant x$ in \eqref{choices} then gives that 
\begin{equation}
\mathcal{S}^y(x) \ll \frac{x}{(\log x)^2}\label{errorbound}.
\end{equation}
\begin{note4.2}
\textnormal{The result holds even for $\ell=0$ and in this case there is only one factorization $n=ab$. The fact that $\ell\neq 0$ is crucial, however, for what follows in \S 4.3.}
\end{note4.2}
\subsection{Treating $\mathcal{S}_{y}(x)$ by the Large Sieve}
From our definition (\ref{Sy}) of $\mathcal{S}_{y}(x)$, we write $n$ and $n+\ell$ in terms of their factorizations $\mathcal{F}_z$ and are left with evaluating
\begin{equation}
\mathcal{S}_{y}(x)=\sum_{\substack{
     a, a_\ell \leqslant y\\
     p|aa_\ell \Rightarrow p\leqslant z
   }} | \lambda_1 ( a ) \lambda_2 ( a_\ell ) |
   \sum_{\substack{
     n \leqslant x \\
     n \equiv 0 \, (\textnormal{mod } a)\\
     n \equiv -\ell \, (\textnormal{mod } a_\ell)\\
     p|bb_\ell \Rightarrow p>z
   }} |\lambda_1 ( b ) \lambda_2(b_\ell)|. \label{writeasfactors}
\end{equation}
To help deal with certain co-primality conditions which come up during analysis, we pull out the greatest common divisor $v$ of $a$ and $a_\ell$, writing now $n_v:=n/v=ab$ and $n_v+w=a_\ell b_\ell$ with $(a,a_\ell)=(aa_\ell,w)=1$ so that 
\begin{equation}
\mathcal{S}_{y}(x) = \sum_{vw=\ell} \sum_{\substack{
     a, a_\ell \leqslant y/v\\
     p|aa_{\ell} \Rightarrow p\leqslant z \\
     (a,a_\ell)=(aa_\ell,w)=1
   }} | \lambda_1 ( v a ) \lambda_2 ( v a_\ell) | \sum_{\substack{
     n_v \leqslant x / v\\
     n_v \equiv r \, (\textnormal{mod } aa_\ell)\\
     p|bb_\ell \Rightarrow p>z
   }} |\lambda_1 ( b ) \lambda_2(b_\ell)|.\label{dividebyv}
\end{equation}
Here we applied the Chinese remainder theorem so that the residue class $r$ in the inner-most sum satisfies $r \equiv 0 \, (\textnormal{mod } a)$ and $r \equiv -w \, (\textnormal{mod } a_\ell)$.  Finally, we take advantage of positivity by applying the Ramanujan-Petersson conjecture to have that 
\begin{equation}
|\lambda_1(b)\lambda_2(b_\ell)| \ll (\log x)^{2m\varepsilon} \label{RPbound}
\end{equation}
by our choice of $s$ in \eqref{choices}.  Indeed, we have $|\lambda_1(p^\alpha)|\leqslant \tau_m(p^\alpha) \leqslant 2^{\alpha+m-1}$ for prime $p$ and $b=p_1^{\alpha_1} p_2^{\alpha_2}\ldots p_t^{\alpha_t}$ with $\alpha_1+\alpha_2+\ldots + \alpha_t \leqslant s$.  

We proceed to bound the count
\begin{equation}
\sum_{\substack{
      n_v \leqslant x / v\\
     n_v \equiv r \, (\textnormal{mod } aa_\ell)\\
     p|bb_\ell \Rightarrow p>z
   }} 1. \label{pulloutRP}
\end{equation}  
Writing $n_v=(aa_\ell)m+r$ with $0\leqslant r < aa_\ell$, we note the following equivalences between divisability conditions for primes $p\leqslant z$
\begin{eqnarray}
p \nmid b  & \Longleftrightarrow & p \nmid (a_\ell m + r/a), \label{condb}\\
p \nmid b_\ell & \Longleftrightarrow & p \nmid (am+(r+w)/a_\ell) \label{condbell}.
\end{eqnarray}
For fixed $(a,a_\ell)=(aa_\ell,w)=1$ with $a,a_\ell \leqslant y/v$,  we see that the count in \eqref{pulloutRP} is bounded by $S=|\mathcal{S}(\mathcal{M}, \mathcal{P}, \Omega)|$ where we define the \textit{sifted set} to be
\begin{equation*}
\mathcal{S}(\mathcal{M}, \mathcal{P}, \Omega) = \{m\in \mathcal{M};\, m \, (\textnormal{mod } p) \notin \Omega_p \textnormal{ for all } p \in \mathcal{P} \}
\end{equation*}
with
\begin{eqnarray*}
\mathcal{M} & := & \{ m \in \mathbbm{Z} \mid 0 < m (v a a_\ell) \leqslant x \},\\
\mathcal{P} & := & \{ p \textnormal{ prime } \mid 2<p\leqslant z \}
\end{eqnarray*}
and the set of residue classes to be ``sieved out" $\Omega := \bigcup_{p\in\mathcal{P}} \,\Omega_p$ is given by
\begin{equation*}
 \Omega_p:= \left\{\begin{array}{ll} \{r_1 \, (\textnormal{mod }p)\} & \textnormal{for } p|a \\
 \{r_2 \, (\textnormal{mod }p)\} & \textnormal{for } p|a_\ell\\
 \{r_1 \, (\textnormal{mod }p), r_2 \, (\textnormal{mod }p)\} & \textnormal{for } p\nmid a a_\ell \end{array}\right.
\end{equation*}
where $r_1\equiv-\overline{a_{\ell}} r/a \, (\textnormal{mod }p)$ and $r_2\equiv-\overline{a} (r+w)/a_{\ell}\, (\textnormal{mod }p)$. Here the overline means multiplicative inverse modulo $p$.  Indeed, if $p|a$ (or $p|a_\ell$) we see that the condition \eqref{condbell} (resp. \eqref{condb}) is redundant with $p \nmid w$ and therefore only one residue class is ``sieved out" in these cases. Recall that $(a, a_\ell)=(aa_\ell,w)=1$ so that $p\leqslant z$ can not divide both $a$ and $a_\ell$. Setting $\omega(p)=|\Omega_p|$ for all $p\in \mathcal{P}$ and applying the Large Sieve as stated in ([I-K], Theorem 7.14), we have
\begin{equation*}
S \leqslant \frac{N+Q^2}{H}
\end{equation*}
for any $Q\geqslant 1$, where $N=x/(vaa_\ell)$,
\begin{equation*}
H=\sum_{q\leqslant Q} h(q) = \sum_{\substack{q_1\leqslant Q \\ q_1|a a_\ell}} h(q_1) \sum_{\substack{q_2 \leqslant Q/q_1 \\ (q_2,a a_\ell)=1}} h(q_2)
\end{equation*}
and $h(q)$ is the multiplicative function supported on square-free integers with prime divisors in $\mathcal{P}$ such that
\begin{equation*}
h(p)=\frac{\omega(p)}{p-\omega(p)}.
\end{equation*}
Choosing say $Q=x^{1/4}$, we get by our values of $\omega(p)$ that
\begin{equation*}
H \gg \left(\frac{\varphi(a a_\ell)}{a a_\ell}\right)^2 \sum_{\substack{q_1\leqslant Q \\ q_1|a a_\ell}} h(q_1) \sum_{q_2 \leqslant Q/q_1} \frac{\tau(q_2)}{q_2} \gg \frac{\varphi(a a_\ell)}{a a_\ell} (\log z)^2
\end{equation*}
so that $S$ and the count in \eqref{pulloutRP} is
\begin{equation*}
\ll \frac{x}{v\varphi(aa_\ell) (\log z)^2}.
\end{equation*}
Plugging the above bound for \eqref{pulloutRP} along with \eqref{RPbound} back into \eqref{dividebyv} gives 
\begin{equation*}
\mathcal{S}_y(x) \ll  \frac{x}{(\log x)^{2-\varepsilon}} \sum_{vw=\ell} \frac{1}{v}\sum_{\substack{
     a, a_\ell \leqslant y/v\\
     p|aa_{\ell} \Rightarrow p\leqslant z \\
     (a,a_\ell)=(aa_\ell,w)=1
   }} \frac{| \lambda_1 ( v a ) \lambda_2 ( v a_\ell) |}{\varphi(aa_\ell)}
\end{equation*}
and therefore,
\begin{equation*}
\mathcal{S}_y(x) \ll  \frac{x}{(\log x)^{2-\varepsilon}} \prod_{p\leqslant z} \left( 1+ \frac{|\lambda_1(p)|}{p}\right)\left( 1+ \frac{|\lambda_2(p)|}{p}\right) \tau(|\ell|). 
\end{equation*}
By our partition of $S_\ell(x)$ in \eqref{2parts}, combining the above bound for $\mathcal{S}_y(x)$ with the bound for $\mathcal{S}^y(x)$ in \eqref{errorbound} concludes the proof of Theorem 2.

\end{document}